\documentstyle[12pt]{article}

\font\twlgot =eufm10 scaled \magstep1
\font\egtgot =eufm8
\font\sevgot =eufm7
\font\twlmsb =msbm10 scaled \magstep1
\font\egtmsb =msbm8
\font\sevmsb =msbm7

\newfam\gotfam

\textfont\gotfam\twlgot
\scriptfont\gotfam\egtgot
\scriptscriptfont\gotfam\sevgot

\newfam\msbfam
\textfont\msbfam\twlmsb
\scriptfont\msbfam\egtmsb
\scriptscriptfont\msbfam\sevmsb
\def\Bbb{\protect\pBbb}
\def\pBbb{\relax\ifmmode\expandafter\Bb\else\typeout{You cann't use
Bbb in text mode}\fi}
\def\Bb #1{{\fam\msbfam\relax#1}}

\def\thebibliography#1{\section*{References}\list
  {[\arabic{enumi}]}{\settowidth\labelwidth{#1}\leftmargin\labelwidth
    \advance\leftmargin\labelsep
    \usecounter{enumi}}
    \def\newblock{\hskip .11em plus .33em minus .07em}
    \sloppy\clubpenalty4000\widowpenalty4000
    \sfcode`\.=1000\relax}

\def\op#1{\mathop{\fam0 #1}\limits}

\newcommand{\id}{{\rm Id\,}}

\newcommand{\beq}{\begin{equation}}
\newcommand{\eeq}{\end{equation}}
\newcommand{\ben}{\begin{eqnarray}}
\newcommand{\een}{\end{eqnarray}}
\newcommand{\be}{\begin{eqnarray*}}
\newcommand{\ee}{\end{eqnarray*}}
\newcommand{\bea}{\begin{eqalph}}
\newcommand{\eea}{\end{eqalph}}
\newcommand{\cA}{{\cal A}}

\newcommand{\cD}{{\cal D}}

\newcommand{\cH}{{\cal H}}
\newcommand{\cF}{{\cal F}}

\newcommand{\bL}{{\bf L}}

\newcommand{\bC}{{\bf C}}

\newcommand{\al}{\alpha}

\newcommand{\bt}{\beta}
\newcommand{\dl}{\delta}
\newcommand{\la}{\lambda}
\newcommand{\La}{\Lambda}
\newcommand{\f}{\phi}

\newcommand{\Om}{\Omega}
\newcommand{\m}{\mu}

\newcommand{\g}{\gamma}
\newcommand{\G}{\Gamma}

\newcommand{\vt}{\vartheta}
\newcommand{\vf}{\varphi}

\newcommand{\lng}{\langle}
\newcommand{\rng}{\rangle}

\newcommand{\si}{\sigma}
\newcommand{\Si}{\Sigma}
\newcommand{\w}{\wedge}
\newcommand{\wt}{\widetilde}
\newcommand{\wh}{\widehat}
\newcommand{\ol}{\overline}
\newcommand{\dr}{\partial}

\newcommand{\ot}{\otimes}

\newcommand{\ve}{\varepsilon}

\newcounter{eqalph}
\newcounter{equationa}
\newcounter{theorem}
\newcounter{remark}
\newcounter{proposition}
\newcounter{lemma}
\newcounter{corollary}
\newcounter{definition}

\newenvironment{eqalph}{\stepcounter{equation}
\setcounter{equationa}{\value{equation}}
\setcounter{equation}{0}

\begin{eqnarray}}{\end{eqnarray}\setcounter{equation}{\value{equationa}}}

\def\theremark{\arabic{remark}}

\def\thedefinition{\arabic{definition}}

\newenvironment{proof}{\noindent 
{\it Proof.}}{\medskip}

\newenvironment{prop}{\refstepcounter{definition} 
\bigskip\noindent{\it Proposition \thedefinition.}}{\medskip}

\newcommand{\mar}[1]{}

\begin{document}
\hbox{}

{\parindent=0pt

{\large\bf Holonomy control operators in classical and
quantum completely integrable Hamiltonian systems}
\bigskip

{\sc Giovanni Giachetta\footnote{E-mail address:
giovanni.giachetta@unicam.it}, Luigi
Mangiarotti$\dagger$\footnote{E-mail 
address: luigi.mangiarotti@unicam.it} and Gennadi 
Sardanashvily$\ddagger$\footnote{E-mail address: sard@grav.phys.msu.su; 
URL: http://webcenter.ru/$\sim$sardan/ }}
\bigskip

\begin{small}
$\dagger$ Department of Mathematics and Physics, University of Camerino, 62032
Camerino (MC), Italy \\
$\ddagger$ Department of Theoretical Physics, Physics Faculty, Moscow State
University, 117234 Moscow, Russia
\bigskip
\bigskip

{\bf Abstract.}

Completely integrable Hamiltonian systems look rather promising
for controllability since their first integrals are stable
under an internal evolution, and
one may hope to find a perturbation of a Hamiltonian
which drives the first integrals
at will. Action-angle coordinates are most convenient
for this purpose. Written with respect to these coordinates,
a Hamiltonian and first integrals of a (time-dependent)
completely integrable system 
depend only on action variables.
We introduce a suitable perturbation of an internal Hamiltonian by a term
containing time-dependent parameters (control fields) so that 
an evolution 
of action variables is nothing else than a holonomy displacement along a curve 
in a parameter space. Therefore, one can determine this evolution
in full by an appropriate choice of control fields.
Similar holonomy controllability of finite level
quantum systems is of special interest in connection with
quantum computation. We provide geometric quantization of
a time-dependent completely integrable Hamiltonian system
in action-angle variables. Its Hamiltonian and first
integrals have time-independent countable spectra. 
The holonomy control operator in this countable level
quantum system is constructed.
\end{small}
}


\section{Introduction}

A time-dependent Hamiltonian system of $m$ degrees of freedom is called
a completely integrable system (henceforth CIS)
if there exist $m$ independent first integrals 
in involution. We show that such a system admits the 
action-angle coordinates around any regular instantly compact invariant manifold.
Written relative to these coordinates, its Hamiltonian and first integrals
are functions only of action coordinates. 
In comparison with perturbations in the  
KAM theorem, we consider perturbations of a Hamiltonian of a CIS
by the term which contains time-dependent parameters.
A generic Hamiltonian of a mechanical system with time-dependent
parameters includes a term which is linear in momenta and the temporal derivatives
of parameter functions \cite{book98,jmp02,jmp00b}.
Then the corresponding evolution operator of action variables 
depends only on a trajectory of the parameter
functions in a parameter space. Therefore, it is determined 
completely by an appropriate choice of these parameter functions 
(control fields), and plays the role of a holonomy control operator.

At present, holonomy control 
operators in quantum systems attract special attention in connection with
quantum computation, based on the generalization
of Berry's phase by means of 
driving a finite level 
degenerate eigenstate of a Hamiltonian over the parameter
manifold \cite{fuj,pach,zan}. Information is encoded in
this degenerate state. Bearing in mind this application,
we quantize a time-dependent CIS in question,
and construct the quantum holonomy control operator
in this system. 

An essential simplification is that, given action-angle
variables, a time-dependent CIS can be quantized as an autonomous one
because its Hamiltonian and first integrals 
relative to these coordinates are 
time-independent. Of course,
the choice of action-angle coordinates by no
means is unique. They are subject to
canonical transformations. Therefore, we employ the
geometric quantization technique \cite{eche98,sni,wood},
which remains equivalent 
under canonical transformations, but
essentially depends on a choice of polarization \cite{blat,raw}.

Geometric quantization of
an autonomous CIS has been studied 
with respect to polarization spanned by Hamiltonian vector fields of
first integrals \cite{myk}. The well-known Simms quantization of
the harmonic oscillator is also of this type.
The problem is that the associated quantum algebra
includes functions which are ill defined on the whole momentum
phase space,
and elements of the carrier space fail to be smooth
sections of the quantum bundle.  Indeed, written with respect to the
action-angle variables, this quantum
algebra consists of functions which are affine in angle coordinates.

We choose a different polarization
spanned by almost-Hamiltonian vector fields of angle variables.
The associated quantum algebra $\cA$ 
consists of smooth functions which are affine in action variables. 
This quantization is equivalent to geometric quantization
of the cotangent bundle 
of the torus $T^m$ with respect to the  
vertical polarization. 
This polarization is known to lead to 
Schr\"odinger quantization. 
We show that $\cA$ possesses a set of
inequivalent representations 
in the separable pre-Hilbert space $\Bbb C^\infty(T^m)$
of smooth complex functions
on $T^m$. In particular, the action operators read
\mar{ci7}\beq
\wh I_k=-i\dr_k -\la_k,  \label{ci7}
\eeq
where $\la_k$ are real numbers which specify different representations of $\cA$.
By virtue of the multidimensional Fourier theorem, 
an orthonormal basis for $\Bbb C^\infty(T^m)$ consists of
functions
\mar{ci15}\beq
\psi_{(n_r)}=\exp[i(n_r\f^r)], \qquad (n_r)=(n_1,\ldots,n_m)\in\Bbb Z^m, 
\label{ci15}
\eeq
where $\f^i$ are cyclic coordinates on $T^m$.
With respect to this basis, the action operators (\ref{ci7}) are 
written as countable diagonal matrices
\mar{ci9}\beq
\wh I_k\psi_{(n_r)}=(n_k-\la_k)\psi_{(n_r)}. \label{ci9}
\eeq
Given the representation (\ref{ci7}),
any polynomial Hamiltonian $\cH(I_k)$ of
a CIS is uniquely quantized as a Hermitian element
$\wh\cH(I_k)=\cH(\wh I_k)$ of the enveloping algebra
 of $\cA$. It has the time-independent countable spectrum
\mar{ww10}\beq
\wh \cH(I_k)\psi_{(n_r)}=E_{(n_r)}\psi_{(n_r)}, \qquad 
E_{(n_r)}=\cH(n_k-\la_k), \qquad n_k \in(n_r). \label{ww10}
\eeq
Since $\wh I_k$ are diagonal, one can also quantize Hamiltonians 
$\cH(I_j)$ which are analytic functions on $\Bbb R^m$.

Quantum CISs look rather 
promising for quantum 
computation. Its Hamiltonian depends only on 
action variables and possesses a
time-independent countable spectrum. 
Moreover, it can be 
made degenerate at will by appropriate canonical transformations. 
We construct the holonomy control operator in this 
countable level quantum system.
Note that one used to study controllability of finite level quantum
systems \cite{alb,alt,sch}.

\section{Classical time-dependent completely integrable systems}

In order to introduce action-angle variables for a time-dependent CIS,
we use the fact that a time-dependent Hamiltonian system of
$m$ degrees of freedom 
can be extended to an 
autonomous one of $m+1$ degrees of
freedom where the time is treated as a dynamic variable \cite{bouq,dew,lich}.
The 
classical theorem \cite{arn,laz} on 
action-angle coordinates around a regular compact invariant manifold
can not be applied to this autonomous CIS since its
invariant manifolds 
are never compact because of the time axis.
Generalizing the above mentioned theorem, we first prove that 
there is a system of action-angle coordinates on an open neighbourhood
$U$ of a regular connected
invariant manifold $M$ of an autonomous CIS
if Hamiltonian vector fields of first integrals on $U$ are complete
and the foliation of $U$ by invariant manifolds is trivial.
If $M$ is a compact regular invariant manifold, these conditions always hold
\cite{laz}. Afterwards, we show that, if a regular connected invariant manifold
of a time-dependent CIS 
is compact at each instant, it
is diffeomorphic to the product 
of the time 
axis $\Bbb R$ and an $m$-dimensional torus $T^m$, and admits
an open neighbourhood endowed with
time-dependent action-angle coordinates $(I_i;t,\f^i)$,
$i=1,\ldots,m$, where $t$ is the Cartesian coordinate on $\Bbb R$ and
$\f^i$ are cyclic coordinates on $T^m$.

Recall that the configuration space of a time-dependent mechanical system
is a fibre bundle $Q\to \Bbb R$
over the time axis $\Bbb R$ equipped with the bundle
coordinates $(t,q^k)$, $k=1,\ldots,m$.
The corresponding momentum phase space is the vertical
cotangent bundle
$V^*Q$ of $Q\to\Bbb R$ endowed with holonomic
coordinates $(t,q^k,p_k=\dot q_k)$ \cite{book98,sard98}.
The cotangent bundle
$T^*Q$, coordinated by $(q^\la,p_\la)=(t,q^k,p_0,p_k)$, 
plays a role of the homogeneous momentum phase space. It is provided 
with the canonical Liouville form $\Xi=p_\la dq^\la$, the
canonical symplectic form
$\Om=dp_\la\w dq^\la$, 
and the corresponding Poisson bracket
\mar{z7}\beq
\{f,f'\}_T =\dr^\la f\dr_\la f' -\dr_\la
f\dr^\la f', \qquad f,f'\in C^\infty(T^*Q). \label{z7}
\eeq
There is the one-dimensional trivial affine bundle
\mar{z11}\beq
\zeta:T^*Q\to V^*Q. \label{z11}
\eeq
Given its global section $h$, one can equip $T^*Q$ with 
the global fibre coordinate $r=p_0-h$.
The fibre bundle (\ref{z11}) provides
the vertical cotangent bundle $V^*Q$ with the canonical Poisson
structure $\{,\}_V$ such that
\mar{m72',2}\ben
&& \zeta^*\{f,f'\}_V=\{\zeta^*f,\zeta^*f'\}_T, \qquad \forall 
f,f'\in C^\infty(V^*Q), \label{m72'} \\
&& \{f,f'\}_V = \dr^kf\dr_kf'-\dr_kf\dr^kf'. \label{m72}
\een

A Hamiltonian of time-dependent mechanics is defined
as a global section
\mar{qqq}\beq
h:V^*Q\to T^*Q, \qquad p_0\circ
h=-\cH(t,q^j,p_j), \label{qqq}
\eeq
of the affine bundle $\zeta$ (\ref{z11}) \cite{book98,sard98}. It
yields the pull-back Hamiltonian form
\be
H=h^*\Xi= p_k dq^k -\cH dt 
\ee
on $V^*Q$. There exists a unique
vector field $\g_H$ on $V^*Q$ such that 
\mar{z3}\ben
&& \g_H\rfloor dt=1, \qquad \g_H\rfloor dH=0, \nonumber \\
&& \g_H=\dr_t + \dr^k\cH\dr_k- \dr_k\cH\dr^k. \label{z3}
\een
Its trajectories obey the Hamilton equation
\mar{z20}\beq
\dot q^k=\dr^k\cH, \qquad \dot p_k=-\dr_k\cH. \label{z20}
\eeq

A first integral of
the Hamilton equation (\ref{z20}) is defined as
a smooth real function $F$ 
on $V^*Q$ whose Lie derivative 
\be
\bL_{\g_H} F=\g_H\rfloor dF=\dr_tF +\{\cH,F\}_V 
\ee
along the vector field $\g_H$ (\ref{z3}) vanishes, i.e., $F$
is constant on trajectories of $\g_H$. 
A time-dependent Hamiltonian system $(V^*Q,H)$ is said to be
completely integrable
if the Hamilton equation (\ref{z20}) admits $m$ first integrals 
$F_k$ which are
in involution with respect to the Poisson bracket $\{,\}_V$ (\ref{m72}),
and whose differentials $dF_k$ are linearly independent almost everywhere
(i.e., the set of points where this condition fails is nowhere dense).
One can associate to this CIS an autonomous CIS on $T^*Q$ as follows.

Given a Hamiltonian $h$ (\ref{qqq}), it is readily observed that 
\mar{mm16}\beq
\cH^*=\dr_t\rfloor(\Xi-\zeta^* h^*\Xi)=p_0+\cH \label{mm16}
\eeq
is a function on $T^*Q$.  
Let us regard $\cH^*$
as a Hamiltonian of an autonomous Hamiltonian system on the symplectic
manifold $(T^*Q,\Om)$ \cite{jmp00}. Its Hamiltonian vector field 
\mar{z5}\beq
\g_T=\dr_t -\dr_t\cH\dr^0+ \dr^k\cH\dr_k- \dr_k\cH\dr^k \label{z5}
\eeq
is projected onto the vector field $\g_H$ (\ref{z3}) on $V^*Q$ so that
\be
\zeta^*(\bL_{\g_H}f)=\{\cH^*,\zeta^*f\}_T, \qquad
\forall f\in C^\infty(V^*Q).
\ee
An immediate consequence of this relation is the following.

\begin{prop} \label{z6} \mar{z6} 
(i) Given a time-dependent CIS $(H; F_k)$ on $V^*Q$, the 
Hamiltonian system
$(\cH^*,\zeta^*F_k)$ on $T^*Q$ is a CIS.
(ii) Let $N$ be a connected regular invariant manifold of $(H; F_k)$. 
Then $h(N)\subset
T^*Q$ is a connected regular invariant manifold of 
the autonomous CIS $(\cH^*,\zeta^*F_k)$.
\end{prop}

Hereafter, we assume that the  
vector field $\g_H$ (\ref{z3}) is complete. In this case, the Hamilton equation
(\ref{z20}) admits a unique global solution through each point of the
momentum phase space $V^*Q$, and trajectories of $\g_H$
define a trivial bundle $V^*Q\to V^*_tQ$
over any fibre 
$V^*_tQ$ of $V^*Q\to \Bbb R$. Without loss of generality, we choose
the fibre $i_0:V^*_0Q\to V^*Q$
at $t=0$.  Since $N$ is an
invariant manifold, the fibration 
\mar{ww}\beq
\xi:V^*Q\to V^*_0Q \label{ww}
\eeq
also yields the fibration of $N$ onto 
$N_0=N\cap V^*_0Q$
such that $N\cong \Bbb R\times N_0$
is a trivial bundle. 

\section{Time-dependent action-angle coordinates}

Let us introduce the  
action-angle coordinates around an
invariant manifold $N$ of a time-dependent CIS on $V^*Q$ by use of the action-angle
coordinates around the invariant manifold $h(N)$ of the
autonomous CIS on $T^*Q$ in
Proposition \ref{z6}. Since $N$ and, consequently, $h(N)$ are non-compact, we 
refer to the following.

\begin{prop} \label{z8} \mar{z8}
Let $M$ be a connected invariant manifold of an autonomous CIS
 $\{F_\la\}$, $\la=1,\ldots,n$, on a symplectic manifold
$(Z,\Om_Z)$. Let $U$ be an open neighbourhood of $M$ such that: (i) the differentials
$dF_\la$ are independent everywhere 
in $U$, (ii) the Hamiltonian vector fields $\vt_\la$ of 
the first integrals $F_\la$ on $U$ are complete,
and (iii) the submersion $\times F_\la: U\to \Bbb R^n$
is a trivial bundle of invariant manifolds 
over a domain $V'\subset \Bbb R^n$. 
Then $U$ is isomorphic
to the symplectic annulus 
\mar{z10}\beq
W'=V'\times(\Bbb R^{n-m}\times T^m), \label{z10}
\eeq
provided with the action-angle coordinates 
\mar{z11'}\beq
(I_1,\ldots,I_n; x^1,\ldots, x^{n-m}; \f^1,\ldots,\f^m) \label{z11'}
\eeq
such that the symplectic form on $W'$ reads
\be
\Om_Z=dI_a\w dx^a +dI_i\w d\f^i,
\ee
and the first integrals $F_\la$ depend only on  
the action coordinates $I_\al$.
\end{prop}

\begin{proof}
In accordance with the well-known theorem \cite{arn}, 
the invariant
manifold $M$ is diffeomorphic to the product $\Bbb R^{n-m}\times T^m$,
which is the group space of the quotient $G=\Bbb R^n/\Bbb Z^m$ 
of the group $\Bbb R^n$ generated by Hamiltonian vector fields
$\vt_\la$ of first integrals $F_\la$ on $M$. 
Namely, $M$
is provided with the group space coordinates $(y^\la)=(s^a,\vf^i)$ 
where $\vf^i$ are linear functions 
of parameters $s^\la$ along integral curves of the 
Hamiltonian vector fields $\vt_\la$ on $U$. Let
$(J_\la)$ be coordinates on
$V'$ which 
are values of first integrals $F_\la$. 
Let us choose a trivialization
of the fibre bundle $U\to V$ seen as a principal bundle with the
structure group $G$. We fix its global section $\chi$.
Since parameters $s^\la$ are given up to a shift, let us
provide each fibre $M_J$, $J\in V$, with the group space
coordinates $(y^\la)$ 
centred at the point 
$\chi(J)$. Then 
$(J_\la;y^\la)$ are bundle coordinates
on the annulus $W'$ (\ref{z10}). The rest of the proof 
in 
Appendix A reduces 
to transformation of the coordinates $(J_\la;y^\la)$ to the desired
coordinates
\mar{ww21}\beq
I_a=J_a, \qquad I_i(J_j), \qquad x^a=s^a + S^a(J_\la), \qquad 
\f^i=\vf^i+ S^i(J_\la,s^b).
\label{ww21}
\eeq 
\end{proof}

Of course, the action-angle coordinates (\ref{z11'}) by no 
means are unique. 
For instance, let $\cF_a$ be $(n-m)$ arbitrary smooth functions
on $\Bbb R^m$. 
Let us consider the canonical coordinate transformation 
\mar{ww26}\beq
I'_a=I_a-\cF_a(I_j), \qquad I'_k=I_k, \qquad x'^a=x^a, \qquad \f'^i= \f^i + 
x^a\dr^i\cF_a(I_j). \label{ww26}
\eeq
Then $(I'_a,I'_k; x'^a,\f'^k)$ are action-angle
coordinates on the symplectic annulus which differs from $W'$ (\ref{z10})
in a trivialization. 

Let us apply Proposition \ref{z8} to the CISs in Proposition \ref{z6}. 

\begin{prop} \label{z13} \mar{z13}
Let $N$ be a connected regular invariant manifold of a time-dependent 
CIS on $V^*Q$, and let the image $N_0$ of its 
projection $\xi$ (\ref{ww}) be compact.
Then the invariant manifold $h(N)$ of the associated autonomous CIS
on $T^*Q$ has an open
neighbourhood $U$ obeying the condition of Proposition \ref{z8} (see Appendix B).
\end{prop}

In accordance with Proposition \ref{z8}, the open neighbourhood $U$ of
the invariant manifold $h(N)$ of the autonomous CIS
in Proposition \ref{z13}
is isomorphic to the symplectic annulus
\mar{z41}\beq
W'=V'\times(\Bbb R\times T^m) \label{z41}
\eeq
provided with the  action-angle coordinates 
$(I_0,\ldots,I_m;t,\f^1,\ldots,\f^m)$
such that 
the symplectic form on 
$W'$ reads
\be
\Om=dI_0\w dt + dI_k\w d\f^k. 
\ee
By the construction in Proposition \ref{z8}, 
$I_0=J_0=\cH^*$ (\ref{mm16}) and the corresponding angle coordinate 
is $x^0=t$, while the first integrals $J_k=\zeta^*F_k$
depend only on the action coordinates $I_i$.

Since the action coordinates $I_i$ are independent of the coordinate
$J_0$, the symplectic annulus $W'$ (\ref{z41}) inherits the fibration
(\ref{z11}) which reads
\mar{z46}\beq
\zeta: W'\ni (I_0,I_i;t,\f^i)\to 
(I_i,t,\f^i)\in W=\Bbb R\times T^m\times V.  \label{z46}
\eeq
By the relation similar to (\ref{m72'}), the product $W$
(\ref{z46}) is provided with the Poisson structure
\be
\{f,f'\}_W = \dr^if\dr_if'-\dr_if\dr^if', \qquad f,f'\in C^\infty(W).
\label{ww2}
\ee
Therefore, one can regard $W$  as the momentum phase space of the 
time-dependent CIS
in question around an invariant manifold $N$. 

It is readily observed that the Hamiltonian vector field $\g_T$ of the 
autonomous Hamiltonian
$\cH^*=I_0$ is $\g_T=\dr_t$, and so is its projection $\g_H$ (\ref{z3})
on $W$. Consequently, the Hamilton equation (\ref{z20})
of a time-dependent CIS
with respect to the action-angle coordinates take the form
$\dot I_i=0$, $\dot\f^i=0$.
Hence, $(I_i;t,\f^i)$ are the initial data coordinates.
One can introduce such coordinates as follows. Given the fibration $\xi$ (\ref{ww}),
let us provide $N_0\times V\subset V^*_0Q$ in Proposition \ref{z13} with 
the action-angle coordinates $(\ol I_i;\ol \f^i)$ for the 
CIS $\{i_0^*F_k\}$ on the symplectic leaf
$V^*_0Q$. Then, it is readily observed that $(\ol I_i;t,\ol \f^i)$ are 
time-dependent action-angle coordinates
on $W$ (\ref{z46}) such that the Hamiltonian 
$\cH(\ol I_j)$ of a time-dependent CIS relative to these coordinates vanishes,
i.e., $\cH^*=\ol I_0$. Using the canonical transformations (\ref{ww26}),
one can obtain different time-dependent action-angle coordinates. 
In particular, given a smooth function $\cH$ on $\Bbb R^m$,
one can endow  $W$ with
the action-angle coordinates 
\be
I_0=\ol I_0-\cH(\ol I_j), \qquad I_i=\ol I_i, \qquad \f^i=\ol \f^i + 
t\dr^i\cH(\ol I_j)
\ee
such that $\cH(I_i)$ is a Hamiltonian of time-dependent CIS on 
$W$. 

\section{The classical control operator}

A generic momentum phase space of a Hamiltonian system
with time-dependent parameters is a composite fibre bundle
$\Pi\to\Si\to\Bbb R$,
where $\Pi\to\Si$ is a symplectic bundle and $\Si\to\Bbb R$ is a parameter
bundle whose sections are parameter functions 
\cite{jmp02,book98,jmp00b,wu}. Here, we assume that all bundles are
trivial and, moreover, their trivializations hold fixed. Then the
momentum phase space of a Hamiltonian system with time-dependent
parameters on
the Poisson manifold $W$ (\ref{z46}) is the product
\mar{zz2}\beq
\Pi= S\times W=V\times (\Bbb R\times S\times T^m)\to \Bbb R\times S\to\Bbb R,
\label{zz2} 
\eeq
equipped with the coordinates $(I_k;t,\si^\al,\f^k)$. It is convenient
to suppose for a time that parameters are also dynamic variables. The
momentum phase space of such a system is
$\Pi'=T^*S\times W$, 
coordinated by $(I_k;t,\si^\al,p_\al,\f^k)$.
The dynamics of a time-dependent mechanical system on the momentum phase
space $\Pi'$ is characterized by a Hamiltonian form
\mar{pr30}\ben
&& H_\Si= p_\al d\si^\al+ I_k d\f^k -
\cH_\Si(t,\si^\bt,p_\bt, I_j,\f^j)dt
\nonumber\\
&& \cH_\Si=p_\al\G^\al_t +I_k(\La^k_t +\G^\al_t\La^k_\al) 
+\wt\cH, \label{pr30}
\een
where $\G=(\G^\al_t)$ is a connection 
on the parameter bundle $\Bbb R\times S\to\Bbb R$ and 
$\La=(\La^k_t,\La^k_\al)$ is a connection 
on $\Bbb R\times S\times T^m\to \Bbb R\times S$ \cite{jmp02,book98,jmp00b}. 

Bearing in mind that $\si^\al$ are parameters, one should choose the
Hamiltonian $\cH_\Si$ (\ref{pr30}) to be affine in momenta $p_\al$.  
Furthermore, in order to describe a Hamiltonian system with a fixed
parameter function  $\si^\al=\xi^\al(t)$, one defines the
connection
$\G$ such that
\be
\nabla^\G \xi=0,\qquad \G^\al_t(t,\xi^\bt(t))=\dr_t\xi^\al. 
\ee
Then the pull-back
\be
H_\xi=\xi^*H_\Si=I_k d\f^k - (I_k[\La^k_t(t,\xi^\bt,\f^j)+ \La^k_\al
(t,\xi^\bt,\f^j)\dr_t\xi^\al]
+\wt\cH(t,\xi^\bt,I_j,\f^j))dt 
\ee
is a Hamiltonian form on the Poisson manifold $W$ (\ref{z46}). 
Let us put
\be
\wt\cH=\cH-I_k\La^k_t,
\ee
where $\cH(I_i)$ is a Hamiltonian of the original CIS
 on $W$ (\ref{z46}).
Then the Hamiltonian form
\mar{zz5}\beq
H_\xi=I_k d\f^k-\cH_\xi dt=I_k d\f^k - [I_k\La^k_\al
(t,\xi^\bt,\f^j)\dr_t\xi^\al+ \cH(I_j)]dt \label{zz5}
\eeq
describes a perturbed time-dependent CIS
on the Poisson manifold $W$ (\ref{z46}). 
The corresponding Hamilton equation reads
\mar{zz6}\beq
\dr_t I_k=-\dr_k\La^j_\al I_j\dr_t\xi^\al, \qquad \dr_t\f^k=\dr^k\cH
+\La^k_\al \dr_t\xi^\al. \label{zz6}
\eeq

In order to make the term
\mar{zz8}\beq
\Delta=I_k\La^k_\al\dr_t\xi^\al \label{zz8}
\eeq
in the perturbed Hamiltonian $\cH_\xi$ (\ref{zz5}) a control operator, let us
assume that the coefficients $\La^k_\al$ of the connection $\La$
are independent 
of time. Then, in view of the trivialization (\ref{zz2}), its
part $(\La^k_\al)$ can be seen as a connection
on the fibre bundle $S\times T^m\to S$. Let us choose the initial 
data action-angle
variables. Then the internal Hamiltonian
$\cH$ of the original CIS 
vanishes, and the Hamilton equation (\ref{zz6}) takes the form
\mar{zz10,1}\ben
&& \dr_t I_i=-I_k\dr_i\La^k_\al \dr_t\xi^\al, \label{zz10}\\
&& \dr_t\f^i=\La^i_\al \dr_t\xi^\al. \label{zz11}
\een
This is the control equation as follows.

Any smooth complex function on the product $\Bbb R\times T^m$ is
represented by a multidimensional Fourier series of functions
$\psi_{(n_r)}$ (\ref{ci15})
on $T^m$ with coefficients the smooth functions on $\Bbb R$.
Let us rewrite the equation (\ref{zz11}) as the countable system of equations
\be
\dr_t\psi_{(n_r)}=i\psi_{(n_r)} n_i\dr_t\f^i=
i\psi_{(n_r)}n_i\La^i_\al(\xi^\bt,\psi_{(m_r)}) \dr_t\xi^\al, \qquad n_i\in(n_r), 
\ee
for functions $\psi_{(n_r)}$ (\ref{ci15}). Let 
\mar{qqq2}\beq
\La^i_\al=\op\sum_{(m_r)} \La^i_{\al(m_r)}(\xi^\bt)\psi_{(m_r)} \label{qqq2}
\eeq
be the Fourier series for $\La^i_\al$.
Since $\psi_{(n_r)}\psi_{(m_r)}=
\psi_{(n_r+m_r)}$, we obtain a countable system of linear ordinary 
differential equations
\mar{zz12,'}\ben
&& \dr_t\psi_{(n_r)}= \op\sum_{(n_q)} [iM_{\al(n_r)}^{(k_r)}(\xi^\bt)
\dr_t\xi^\al] \psi_{(k_q)}, \label{zz12}\\
&& M_{\al(n_r)}^{(k_r)} = \op\sum_{(m_r+k_r=n_r)} n_i\La^i_{\al(m_r)}, \label{zz12'}
\een
with time-dependent coefficients $[iM_{\al(n_r)}^{(k_r)}(\xi^\bt)
\dr_t\xi^\al]$. 
Its solution with the initial data $\psi^i(0)$ can be written as the
formal time-ordered exponential 
\mar{zz14}\ben
&&\psi_{(n_r)}(t)=U(t)_{(n_r)}^{(k_r)} \psi_{(k_r)}(0), \nonumber\\
&& U(t)=T\exp\left[i\op\int^t_0 \wh M_\al(\xi^\bt(t'))\dr_t\xi^\al dt'\right]
=T\exp\left[i\op\int_{\xi([0,t])}\wh M_\al(\si^\bt)d\si^\al\right], \label{zz14}
\een
where $\wh M_\al$ denotes the matrix with time-dependent entries 
(\ref{zz12'})
\cite{lam,oteo}. A glance at the evolution operator $U(t)$
(\ref{zz14}) shows that solutions of the equations (\ref{zz12}) 
are functions $\psi_{(n_r)}(\xi(t))$ of a point $\xi(t)$ of the curve 
$\xi:\Bbb R\to S$ in the parameter space $S$.

Substituting this solution into the equation (\ref{zz10}), we obtain the
system of $l$ ordinary linear differential equations with time-dependent
coefficients:
\be
\dr_t I_i=-[L_{\al i}^k (\xi^\bt(t),\psi_{(n_r)}(\xi(t)))
 \dr_t\xi^\al]I_k, \qquad L_{\al i}^k=\dr_i\La^k_\al,
\ee
Its solution is given by the time-ordered
exponential 
\be
&& I_i(t)=U(t)_i^k I_k(0), \\
&& U(t)=T\exp\left[-\op\int^t_0 \wh L_\al(\xi^\bt(t),\psi_{(n_r)}(\xi(t)))\dr_t\xi^\al
dt'\right]= \\
&& \qquad  T\exp\left[-\op\int_{\xi([0,t])}\wh L_\al(\si^\bt,\psi_{(n_r)}(\si)
)d\si^\al\right], 
\ee
where $\wh L_\al$ denotes the matrix with time-dependent entities $L_{\al i}^k$.
This solution is a functions of a point $\xi(t)$ of the curve $\xi:\Bbb R\to
S$ in the parameter space $S$. 

It follows that, if the holonomy group of 
the connection $(\La_\al^k)$ on the fibre bundle $S\times T^m\to S$ is the whole
group $GL(m,\Bbb R)$, one can obtain any desired trajectory
of action variables $I_i$ and, consequently, of first integrals $F_i$
by an appropriate choice of control
functions $\xi^\al(t)$. 

\section{Quantum completely integrable systems}

In order to quantize a time-dependent CIS on the
Poisson manifold $(W,\{,\}_W)$, one may follow the general procedure of
instantwise 
geometric quantization of time-dependent Hamiltonian systems in \cite{sard02}.
As was mentioned above, it can however be quantized
as an autonomous CIS on the symplectic annulus 
\be
P=V\times T^m, 
\ee
equipped with fixed action-angle coordinates $(I_i,\f^i)$ and 
provided with the symplectic form 
\mar{ci1}\beq
\Om_P=dI_i\w d\f^i. \label{ci1}
\eeq

In accordance with the standard geometric quantization procedure \cite{sni,wood},
because the symplectic form $\Om_P$ (\ref{ci1}) is exact, the prequantum bundle
is defined as a trivial complex line bundle $C$ over $P$. 
Since the action-angle coordinates are canonical for the symplectic
form (\ref{ci1}), the prequantum bundle $C$ need no metaplectic 
correction.
Let its trivialization
\mar{ci3}\beq
C\cong P \times \Bbb C \label{ci3}
\eeq
hold fixed. Any other trivialization leads to
equivalent quantization of $P$.
Given the associated bundle coordinates $(I_k;\f^k,c)$, $c\in\Bbb C$, on 
$C$ (\ref{ci3}),
one can treat its sections as smooth complex functions on
$P$.

The Konstant--Souriau prequantization formula associates to
each smooth real function $f\in C^\infty(P)$ on
$P$ the first order differential operator
\mar{lqq46}\beq
\wh f=-i\nabla_{\vt_f} + f \label{lqq46}
\eeq
on sections of $C$, where $\vt_f=\dr^kf\dr_k -\dr_kf\dr^k$
is the Hamiltonian vector field of $f$ and
$\nabla$ is the covariant differential with respect to a
suitable $U(1)$-principal connection on $C$. This connection
preserves the
Hermitian metric $g(c,c')=c\ol c'$ on $C$, and
its curvature form obeys the prequantization
condition $R=i\Om_P$. This connection reads
\mar{ci20}\beq
A=A_0 +icI_kd\f^k\ot\dr_c, \label{ci20}
\eeq
where $A_0$ is a flat $U(1)$-principal connection on $C\to
P$. The equivalence classes of flat
principal connections on $C$ are indexed
by the set $\Bbb R^m/\Bbb Z^m$ of homomorphisms of the de Rham cohomology
group
\be
H^1(P)=H^1(T^m)=\Bbb R^m
\ee
of $P$ to the cycle group $U(1)$ \cite{eche98}.
We choose their representatives of the form
\be
A_0[(\la_k)]=dI_k\ot\dr^k + d\f^k\ot(\dr_k +i\la_kc\dr_c),
\qquad \la_k\in [0,1).
\ee
Then the connection (\ref{ci20}) up to gauge transformations 
reads
\mar{ci14}\beq
A[(\la_k)]=dI_k\ot\dr^k + d\f^k\ot(\dr_k +i(I_k+\la_k)c\dr_c).
  \label{ci14}
\eeq
For the sake of simplicity, we will assume that the numbers $\la_k$
in the expression(\ref{ci14}) belong to $\Bbb R$, but will bear in mind that
connections $A[(\la_k)]$ and $A[(\la'_k)]$ with $\la_k-\la'_k\in\Bbb Z$
are gauge conjugated. Given a connection (\ref{ci14}),
the prequantization operators (\ref{lqq46}) read
\mar{ci4}\beq
\wh f=-i\vt_f +(f-(I_k+\la_k)\dr^kf). \label{ci4}
\eeq

Let us choose the above mentioned angle polarization $V\pi$ which is 
the vertical tangent bundle of the fibration $\pi:P\to T^m$, and
is spanned by the vectors $\dr^k$. 
It is readily observed that the corresponding quantum algebra 
$\cA\subset C^\infty(P)$
consists of affine functions
\mar{ci13}\beq
f=a^k(\f^j)I_k + b(\f^j) \label{ci13}
\eeq
of action coordinates $I_k$. 
The carrier space of its representation by operators (\ref{ci4}) is 
defined as the
space $E$ of sections $\rho$ of the prequantum bundle $C$ of 
compact support
which obey the condition $\nabla_\vt\rho=0$ for any Hamiltonian vector field
$\vt$ subordinate to the distribution $V\pi$. This condition reads
\be
\dr_kf\dr^k\rho=0, \qquad \forall f\in C^\infty(T^m).
\ee
It follows that elements of $E$ are independent of action variables and,
consequently, fail to be of compact support, unless $\rho=0$.
This well-known problem of Schr\"odinger geometric quantization 
is solved as follows \cite{blat2,sard02}.

Let us fix a slice $i_T:T^m\to T^m\times V$.
Let $C_T=i^*_TC$ be the pull-back of the prequantum bundle $C$ (\ref{ci3})
over the torus $T^m$. It is a trivial complex line bundle $C_T=T^m\times\Bbb C$
provided with the pull-back Hermitian metric
$g(c,c')=c\ol c'$. Its sections are smooth complex functions on
$T^m$. Let
\be
\ol A = i^*_TA= d\f^k\ot(\dr_k +i(I_k+\la_k)c\dr_c)
\ee
be the pull-back of the connection $A$ (\ref{ci14}) onto $C_T$.
Let $\cD$ be a metalinear bundle
of complex half-forms on the torus $T^m$.
It admits the canonical lift
of any vector field $\tau$ on $T^m$, and
the corresponding
Lie derivative of its sections reads
\be
\bL_\tau=\tau^k\dr_k+\frac12\dr_k\tau^k.
\ee
Let us consider the tensor product
\mar{ci6}\beq
Y=C_T\ot\cD\to T^m. \label{ci6}
\eeq
Since the Hamiltonian vector fields
\be
\vt_f=a^k\dr_k-(I_r\dr_ka^r +\dr_kb)\dr^k
\ee
of functions $f$ (\ref{ci13}) are projectable onto $T^m$, one can
associate to each
element $f$ of the quantum algebra $\cA$ the first order
differential operator
\mar{lmp135}\beq
\wh f=(-i\ol\nabla_{\pi\vt_f} +f)\ot\id+\id\ot\bL_{\pi \vt_f}=
-ia^k\dr_k-\frac{i}{2}\dr_ka^k-a^k\la_k +b \label{lmp135}
\eeq
on sections of $Y$. A direct
computation shows
that the operators (\ref{lmp135}) obey the Dirac condition
\be
[\wh f,\wh f']=-i\wh{\{f,f'\}}. 
\ee
Sections $\rho_T$ of the quantum bundle $Y\to T^m$ (\ref{ci6})
constitute a pre-Hilbert space $E_T$ with respect to the non-degenerate
Hermitian
form
\be
\lng \rho_T|\rho'_T\rng=\left(\frac1{2\pi}\right)^m\op\int_{T^m}
\rho_T \ol \rho'_T, \qquad \rho_T,\rho'_T\in E_T.
\ee
Then it is readily observed that $\wh f$ (\ref{lmp135}) are Hermitian operators
in $E_T$. In particular, the action operators take the form (\ref{ci7}).

Of course, the above quantization depends on the choice of a
connection $A[(\la_k)]$ (\ref{ci14}) and
a metalinear bundle $\cD$. The latter need not be trivial.
If
$\cD$ is trivial, sections of the quantum bundle $Y\to T^m$ (\ref{ci6})
obey the transformation
rule
\be
\rho_T(\f^k+2\pi)=\rho_T(\f^k)
\ee
for all indices $k$. They are naturally complex smooth functions on $T^m$.
In this case, $E_T$ is the above mentioned pre-Hilbert space 
$\Bbb C^\infty(T^m)$ of complex smooth functions on $T^m$ whose basis 
consists of functions (\ref{ci15}). The action operators $\wh I$ (\ref{ci7})
with respect to this basis are represented by
countable diagonal matrices (\ref{ci9}), while functions $a(\f)$
are decomposed
into the pull-back functions $\psi_{(n_r)}$
which act on $\Bbb C^\infty(T^m)$ by multiplications
\mar{ci11}\beq
\psi_{(n_r)} \psi_{(n'_r)}=\psi_{(n_r)} 
\psi_{(n'_r)}=\psi_{(n_r+n'_r)}. \label{ci11}
\eeq

If $\cD$ is a non-trivial metalinear bundle, sections of the quantum bundle
$Y\to T^m$ (\ref{ci6}) obey the transformation
rule
\mar{ci8}\beq
\rho_T(\f^j+2\pi)=-\rho_T(\f^j)  \label{ci8}
\eeq
for some indices $j$. In this case, the orthonormal basis of the 
pre-Hilbert space
$E_T$ can be represented by double-valued complex functions
\mar{ci10}\beq
\psi_{(n_i,n_j)}=\exp[i(n_i\f^i+ (n_j+\frac12)\f^j)] \label{ci10}
\eeq
on $T^m$. They are eigenvectors
\be
\wh I_i\psi_{(n_i,n_j)}=(n_i-\la_i)\psi_{(n_i,n_j)}, \qquad
\wh I_j\psi_{(n_i,n_j)}=(n_j-\la_j +\frac12)\psi_{(n_i,n_j)}
\ee
of the operators $\wh I_k$ (\ref{ci7}), and the functions 
$a(\f)$
act on the basis (\ref{ci10}) by the above law (\ref{ci11}).
It follows that the representation of $\cA$ 
determined by the connection
$A[(\la_k)]$ (\ref{ci14}) in the space of sections
(\ref{ci8}) of a non-trivial quantum bundle $Y$ (\ref{ci6})
is equivalent to its representation determined by the connection
$A[(\la_i,\la_j-\frac12)]$ in the space $\Bbb C^\infty(T^m)$
of smooth
complex functions on $T^m$.

Therefore, one can restrict the study of representations of
the quantum algebra $\cA$ to its representations in $\Bbb C^\infty(T^m)$ 
associated
to different connections (\ref{ci14}). These representations are 
inequivalent, unless $\la_k-\la'_k\in\Bbb Z$ for all indices $k$.

Now, in order to quantize the Poisson manifold $(W,\{,\}_W)$, one can simply 
replace functions on $T^m$ with those on $\Bbb R\times
T^m$ \cite{sard02,sni}. Let us choose the angle polarization of $W$ 
spanned by the vectors $\dr^k$. 
The corresponding quantum algebra $\cA_W\subset C^\infty(W)$ consists of
affine functions
\mar{ww7}\beq
f=a^k(t,\f^j)I_k + b(t,\f^j)  \label{ww7}
\eeq 
of action coordinates $I_k$, represented by the operators (\ref{lmp135})
in the space $\Bbb C^\infty(\Bbb R\times T^m)$ of smooth complex functions on
$\Bbb R\times T^m$. This space is provided with the structure of the
pre-Hilbert $\Bbb C^\infty(\Bbb R)$-module with respect to the non-degenerate 
$\Bbb C^\infty(\Bbb R)$-bilinear form
\be
\lng \psi|\psi'\rng=\left(\frac1{2\pi}\right)^m\op\int_{T^m}  \psi \ol \psi', 
\qquad \psi,\psi'\in \Bbb C^\infty(\Bbb R\times T^m).  
\ee 
Its basis consists of the pull-backs onto
$\Bbb R\times T^m$ of the 
functions $\psi_{(n_r)}$ (\ref{ci15}).

Since the Poisson structure (\ref{ww2}) defines no dynamics 
on the momentum phase space $W$ (\ref{z46}), 
we should 
quantize the homogeneous momentum
phase space $W'$ (\ref{z41}) in order to
describe evolution of a quantum time-dependent CIS. Following the
general scheme in \cite{sard02,jmp02},  one can provide the relevant geometric 
quantization of the symplectic annulus $(W',\Om')$. 
The corresponding quantum algebra $\cA_{W'}\subset C^\infty(W')$ consists of
affine functions
\be
f=a^\la(t,\f^j)I_\la + b(t,\f^j) 
\ee 
of action coordinates $I_\la$. It suffices to consider its subalgebra
consisting of the elements $f$ and $I_0+f$ for all $f\in \cA_W$ (\ref{ww7}).
They are represented by the operators $\wh f$ (\ref{lmp135}) and $I_0=-i\dr_t$
in the pre-Hilbert module $\Bbb C^\infty(\Bbb R\times T^m)$.
 If a Hamiltonian $\cH(I_j)$ 
of the time-dependent CIS is a polynomial
(or analytic) function in action variables, the Hamiltonian $\cH^*$ of the
associated autonomous CIS is quantized as 
\be
\wh\cH^*=-i\dr_t +\cH(\wh I_j).
\ee  
Then we obtain the Schr\"odinger equation
\be
\wh\cH^*\psi=-i\dr_t\psi +\cH(-i\dr_k -\la_k)\psi=0, \qquad 
\psi\in  \bC^\infty(\Bbb R\times T^m). 
\ee
Its solutions are the series
\be
\psi=\op\sum_{(n_r)} B_{(n_r)} \exp[-iE_{(n_r)}t]\psi_{(n_r)}, \qquad 
B_{(n_r)}\in\Bbb C,
\ee
where $E_{(n_r)}$ are the eigenvalues (\ref{ww10})
of the Hamiltonian $\wh\cH$.

\section{The quantum control operator}

In comparison with the classical control operator in Section 4, we will
construct the quantum control operator which preserves the eigenvalues
of a Hamiltonian, and acts in its degenerate eigenspaces.
For instance, let us 
choose action-angle coordinates
such that 
a Hamiltonian
$\cH$ of a CIS equals $I_1$, and it is independent of other action variables 
$I_a$ ($a,b,c=2,\ldots,m$). It is quantized as $\wh\cH=\wh I_1$.
Its eigenvalues are countably
degenerate. Let us consider the perturbed Hamiltonian 
$\cH_\xi$ (\ref{zz5}) where the perturbation term $\Delta$ (\ref{zz8})
depends only on the action-angle 
coordinates with the indices
$a,b,c=2,\ldots,m$, i.e.,
\be
\cH_\xi=I_1+\La^a_\bt(\xi^\m,\f^b)\dr_t\xi^\bt I_a .
\ee
The perturbation term
\be
\Delta=\La^a_\bt(\xi^\m,\f^b)\dr_t\xi^\bt I_a 
\ee
of this Hamiltonian is an element of the quantum algebra $\cA_W$, and is
quantized by the operator
\be
&&\wh\Delta=-i\La^a_\bt\dr_t\xi^\bt\dr_a -\frac{i}{2}\dr_a(\La^a_\bt)
\dr_t\xi^\bt -\la_a \La^a_\bt\dr_t\xi^\bt=\wh\Delta_\bt\dr_t\xi^\bt,\\
&& \wh\Delta^{(n_1,k_a)}_{\bt(n_1,n_a)}=\op\sum_{(m_a +k_a=n_a)}
[(k_a + \frac12 m_a -\la_a)\La^a_{\bt(m_a)}(\xi^\al)], 
\ee
where the Fourier series decomposition (\ref{qqq2}) is used.

Since the operators $\wh\Delta$ and $\wh\cH$ mutually commute, the
corresponding quantum evolution operator reduces to the product
\mar{zz23}\beq
T\exp\left[-i\op\int_0^t\wh\cH_\xi dt'\right]=
U_1(t)\circ U_2(t)=
T\exp\left[-i\op\int_0^t\wh\cH dt'\right]\circ 
T\exp\left[-i\op\int_0^t\wh\Delta dt'\right].
\label{zz23}
\eeq
The first factor in this product is the dynamic evolution operator of
the quantum CIS. It reads
\mar{zz24}\beq
U_1(t)\psi_{(n_1,n_a)}=
\exp[-i(n_1-\la_1)t]\psi_{(n_1,n_a)}. \label{zz24}
\eeq
Its eigenvalues are countably degenerate.
Recall that the operator
(\ref{zz24}) acts in the pre-Hilbert module $\Bbb
C^\infty(\Bbb R\times T^m)$. Its eigenvalues
are smooth complex functions on $\Bbb R$, and its eigenspaces
are $\Bbb C^\infty(\Bbb R)$-submodules of $\Bbb
C^\infty(\Bbb R\times T^m)$ of countable rank.

The second factor in the product (\ref{zz23}) is
\mar{zz25}\beq
U_2(t)=T\exp\left[-i\op\int_0^t\wh\Delta_\bt(\xi^\al(t'))\dr_t\xi^\bt dt'\right]
 =T\exp\left[-i\op\int_{\xi([0,t])}\wh\Delta_\bt(\si^\al) d\si^\bt\right]. \label{zz25}
\eeq
It acts as a matrix of countable rank in the eigenspaces of the 
internal Hamiltonian $\wh\cH$.   
Its eigenspace corresponding to the eigenvalue $(n_j-\la_j)$ is the 
pre-Hilbert $\Bbb C^\infty(\Bbb
R)$-submodule of 
$\Bbb
C^\infty(\Bbb R\times T^m)$ whose orthonormal basis is made up by  
functions $\psi_{(n_1,n_a)}$ for all collections of integers $(n_a)$.

A glance at the expression (\ref{zz25}) shows that, in fact, the operator
$U_2(t)$ depends on the curve 
$\xi([0,1])\subset S$ in the parameter space $S$. One can treat it as
an operator of parallel
displacement with respect to a connection in the $\Bbb C^\infty(\Si)$-module 
of smooth complex functions on $\Si\times T^m$ along the curve
$\xi$ \cite{jmp02,book00,jmp00}.
For instance, if $\xi([0,1])$ is a loop in $S$, the operator $U_2$
(\ref{zz25}) is  
the geometric Berry factor. In this case, one can think of $U_2$ as
being a holonomy control operator. 

It should be emphasized that, in comparison with 
operators usually studied \cite{alb,alt,fuj,pach,sch,zan}, 
the operator (\ref{zz25}) acts in a countable level quantum
system. Of course, the problem arises if such a system admits complete
controllability and if the holonomy control operator (\ref{zz25})
can provide this controllability. This problem will be studied elsewhere.

\section{Appendix A}

Let us complete the proof of \ref{z8}.
Since $M_J$ are Lagrangian manifolds, 
the symplectic form $\Om_Z$ on $W'$
is given relative to the
bundle coordinates $(J_\la;y^\la)$ by the expression
\mar{ac1}\beq
\Om_Z=\Om^{\al\bt}dJ_\al\w dJ_\bt + \Om^\al_\bt dJ_\al\w dy^\bt. \label{ac1}
\eeq
 By the very definition of coordinates $(y^\la)$, the
Hamiltonian vector fields $\vt_\la$ of first integrals take the
coordinate form $\vt_\la=\vt_\la^\al(J_\m)\dr_\al$. Moreover, since
the cyclic group $S^1$ can not act transitively on $\Bbb R$, we have
\mar{ww25}\beq
\vt_a=\dr_a +\vt_a^i(J_\la)\dr_i, \qquad \vt_i=\vt_i^k(J_\la)\dr_k. \label{ww25}
\eeq
The Hamiltonian vector fields $\vt_\la$ obey the relations
\mar{ww22}\beq
\vt_\la\rfloor\Om_Z=-dJ_\la,\qquad
\Om^\al_\bt \vt^\bt_\la=\dl^\al_\la. \label{ww22}
\eeq
It follows that $\Om^\al_\bt$ is a non-degenerate matrix and 
$\vt^\al_\la=(\Om^{-1})^\al_\la$, i.e., the functions $\Om^\al_\bt$
depend only on coordinates $J_\la$.  
A substitution of (\ref{ww25}) into (\ref{ww22}) results in the equalities 
\mar{ww30,1}\ben
&& \Om^a_b=\dl^a_b, \qquad \vt_a^\la\Om^i_\la=0, \label{ww30}\\
&& \vt^k_i\Om^j_k=\dl^j_i, \qquad \vt^k_i\Om^a_k=0. \label{ww31}
\een
The first of the equalities (\ref{ww31}) shows that the matrix $\Om^j_k$
is non-degenerate, and so is the matrix $\vt^k_i$. Then the second 
one gives $\Om^a_k=0$.

By virtue of
the well-known K\"unneth formula for the de Rham cohomology of a product of manifolds,
the closed form $\Om_Z$ (\ref{ac1}) on $W'$ (\ref{z10})
is exact, i.e., $\Om_Z=d\Xi$ where $\Xi$  reads
\be
\Xi=\Xi^\al(J_\la,y^\la)dJ_\al + \Xi_i(J_\la) d\vf^i + 
\dr_\al\Phi(J_\la,y^\la)dy^\al, 
\ee
where $\Phi$ is a function on $W'$.
Taken up to an exact form, $\Xi$ is brought into the form
\mar{ac2}\beq
\Xi=\Xi'^\al(J_\la,y^\la)dJ_\al + \Xi_i(J_\la) d\vf^i. \label{ac2}
\eeq
Owing to the fact that
components of $d\Xi=\Om_Z$ are independent of $y^\la$ and obey the equalities
(\ref{ww30}) -- (\ref{ww31}), we obtain the following.
 
(i) $\Om^a_i=-\dr_i\Xi'^a +\dr^a\Xi_i=0$. It follows that $\dr_i\Xi'^a$ is 
independent of $\vf^i$, i.e., $\Xi'^a$ is affine in $\vf^i$ and, consequently,
is independent of $\vf^i$ since $\vf^i$ are cyclic coordinate. Hence, 
$\dr^a\Xi_i=0$, i.e., $\Xi_i$ is a function only of coordinates $J_j$.

(ii) $\Om^k_i=-\dr_i\Xi'^k +\dr^k\Xi_i$. Similarly to item (i), one shows that $\Xi'^k$
is independent of $\vf^i$ and $\Om^k_i=\dr^k\Xi_i$, i.e., 
$\dr^k\Xi_i$ is a non-degenerate matrix.

(iii) $\Om^a_b=-\dr_b\Xi'^a=\dl^a_b$. Hence, $\Xi'^a=-s^a+D^a(J_\la)$.

(iv) $\Om^i_b=-\dr_b\Xi'^i$, i.e., $\Xi'^i$ is affine in $s^a$.

In view of items (i) -- (iv), the Liouville form $\Xi$ (\ref{ac2}) reads
\be
\Xi=x^adJ_a + [D^i(J_\la) + B^i_a(J_\la)s^a]dJ_i + \Xi_i(J_j) d\vf^i,
\ee
where we put
\mar{ee1}\beq
x^a=-\Xi'^a=s^a-D^a(J_\la). \label{ee1}
\eeq
Since the matrix $\dr^k\Xi_i$ is non-degenerate,
one can introduce new coordinates $I_i=\Xi_i(J_j)$, $I_a=J_a$. Then we have
\be
\Xi=-x^adI_a + [D'^i(I_\la) + B'^i_a(I_\la)s^a]dI_i + I_i d\vf^i.
\ee
Finally, put 
\mar{ee2}\beq
\f^i=\vf^i-[D'^i(I_\la) + B'^i_a(I_\la)s^a] \label{ee2}
\eeq
in order to obtain the desired action-angle coordinates (\ref{ww21}).
These are bundle coordinates on $U\to V'$ where the 
coordinate shifts (\ref{ee1}) -- (\ref{ee2}) correspond to a choice of another

\section{Appendix B}

In order to prove Proposition \ref{z13}, we first show that 
functions $i_0^*F_k$ make up a CIS on the symplectic
leaf $(V^*_0Q,\Om_0)$ and 
$N_0$ is its invariant manifold without critical points
(i.e., where first integrals fail to be independent).
Clearly, 
the functions $i_0^*F_k$ are in involution, and $N_0$ is their
connected invariant manifold. Let us show that
the set of critical points of $\{i_0^*F_k\}$ is nowhere
dense in $V^*_0Q$ and $N_0$ has none of these points.
Let $V^*_0Q$ be equipped with some coordinates $(\ol q^k,\ol p_k)$.
Then the trivial bundle $\xi$ (\ref{ww}) is provided with the bundle
coordinates $(t,\ol q^k,\ol p_k)$ which play a role of
the initial data coordinates on the momentum phase space $V^*Q$.
Written with respect to these coordinates, the first integrals
$F_k$ become time-independent. It follows that 
\mar{ww12}\beq
dF_k(y)=di_0^*F_k(\xi(y)) \label{ww12}
\eeq
for any point
$y\in V^*Q$. In particular, if $y_0\in V^*_0Q$ is a critical point
of $\{i_0^*F_k\}$, then the trajectory $\xi^{-1}(y_0)$ is a critical set
for the first integrals $\{F_k\}$. The desired statement at once
follows from this result.

 Since $N_0$ is compact and regular,
there is an open neighbourhood of $N_0$
in $V^*_0Q$ isomorphic to $V\times N_0$ where $V\subset \Bbb R^m$ is
a domain, and $\{v\}\times N_0$, $v\in V$, are also invariant
manifolds in $V^*_0Q$ \cite{laz}. Then 
\mar{z70}\beq
W''=\xi^{-1}(V\times N_0) \cong V\times N \label{z70}
\eeq
is an open neighbourhood in
$V^*Q$ of the invariant manifold $N$ foliated by
invariant manifolds $\xi^{-1}(\{v\}\times N_0)$,
$v\in V$, of the time-dependent CIS on $V^*Q$. By virtue of 
the equality (\ref{ww12}), the first integrals $\{F_k\}$ have no 
critical points in $W''$.
For any real number $r\in(-\ve,\ve)$, let us consider a section 
\be
h_r:V^*Q\to T^*Q, \qquad p_0\circ
h_r=-\cH(t,q^j,p_j) +r,
\ee
of the affine bundle $\zeta$ (\ref{z11}). Then the images
$h_r(W'')$ of 
$W''$ (\ref{z70}) make up an open neighbourhood $U$ of
$h(N)$ in $T^*Q$. Because $\zeta(U)=W''$, the
pull-backs $\zeta^*F_k$ of first integrals $F_k$ are free from critical points
in $U$, and so is the function $\cH^*$ (\ref{mm16}). Since
the coordinate $r=p_0-h$ 
provides a trivialization of the affine bundle $\zeta$, the open neighbourhood
$U$ of $h(N)$ is diffeomorphic to the product
\be
(-\ve,\ve)\times h(W'')\cong  
(-\ve,\ve)\times V\times h(N)
\ee
which is a trivialization of the fibration
\be
\cH^*\times(\times \zeta^*F_k): U\to (-\ve,\ve)\times V.
\ee

It remains to prove that the Hamiltonian vector fields of $\cH^*$ and
$\zeta^*F_k$ on $U$ are complete. It is readily observed that the
Hamiltonian vector field $\g_T$ (\ref{z5}) of $\cH^*$ is tangent to
the manifolds $h_r(W'')$, and is the image 
$\g_T=Th_r\circ \g_H\circ \zeta$
of the vector field $\g_H$ (\ref{z3}).
The latter is complete on $W''$, and so is $\g_T$ on
$U$. Similarly, the Hamiltonian vector field
\be
\g_k=-\dr_tF_k\dr^0 +\dr^iF_k\dr_i -\dr_iF_k\dr^i
\ee 
of the function $\zeta^*F_k$ on $T^*Q$ with respect to the Poisson bracket
$\{,\}_T$ (\ref{z7}) is tangent 
to the manifolds $h_r(W'')$, and is the image 
$\g_k=Th_r\circ \vt_k\circ \zeta$
of the Hamiltonian vector field $\vt_k$
of the first integral $F_k$ on $W''$ with respect to
the Poisson bracket $\{,\}_V$ (\ref{m72}). The vector fields $\vt_k$
on $W''$ are vertical relative to the fibration
$W''\to\Bbb R$, and are tangent to compact manifolds.
Therefore, they are complete, and so are the vector fields $\g_k$ on $U$.
Thus, $U$ is the desired open neighbourhood 
of the invariant manifold $h(N)$.

\end{document}